\documentclass{article}

\usepackage{arxiv}

\usepackage[utf8]{inputenc} 
\usepackage[T1]{fontenc}    
\usepackage{hyperref}       
\usepackage{url}            
\usepackage{booktabs}       
\usepackage{amsfonts}       
\usepackage{nicefrac}       
\usepackage{microtype}      
\usepackage{lipsum}		
\usepackage{graphicx}
\usepackage{natbib}
\usepackage{doi}

\usepackage{soul}
\usepackage[utf8]{inputenc}
\usepackage{graphicx}
\usepackage{amsmath}
\usepackage{amsthm}
\usepackage{booktabs}
\usepackage{txfonts}
\usepackage{epsfig}

\usepackage{amssymb}
\usepackage{epsfig}
\usepackage{mathtools}
\usepackage[boxruled]{algorithm2e}

\usepackage{multirow}
\usepackage{booktabs}
\usepackage{balance}
\usepackage{subfig}

\newtheorem{theorem}{Theorem}
\newtheorem{definition}{Definition}
\newtheorem{lemma}{Lemma}

\usepackage[colorinlistoftodos]{todonotes}

\title{Causal Homotopy}

\date{} 					

\author{

       Sridhar Mahadevan\\
   Adobe Research, 345 Park Avenue, San Jose, CA 95110 \\
   
   smahadev@adobe.com
}

\begin{document}

\maketitle

\begin{abstract}
We characterize homotopical equivalences between causal DAG models, exploiting the close connections between  partially ordered set representations of DAGs (posets) and finite Alexandroff topologies. Alexandroff spaces yield a directional topological space: the topology is defined by a unique minimal basis defined by open sets $U_x$ for each variable $x$, specified as the intersection of all open sets containing $x$. Alexandroff spaces induce a (reflexive, transitive) preorder: a variable $x \leq y$ if  $x \in U_y$. Alexandroff spaces satisfying the Kolmogorov $T_0$ separation criterion, where open sets distinguish variables, converts the preordering into a partial ordering. Our approach broadly is to construct a topological representation of posets from data, and then use the poset representation to build a conventional DAG-oriented causal model. We illustrate our framework by showing how it unifies disparate algorithms and case studies proposed previously. Topology plays two key roles in causal discovery. First, topological separability constraints on datasets have been used in several previous approaches to infer causal structure from observations and interventions. Second, a diverse range of graphical models used to represent causal structures can be represented in a unified way in terms of a topological representation of the induced poset structure. We show that the homotopy theory of Alexandroff spaces can be exploited to significantly efficiently reduce the number of possible DAG structures, reducing the search space by several orders of magnitude. 
\end{abstract}


\section{Introduction}

Topology \citep{munkres:book} has found extensive use in many areas in AI, machine learning and optimization.  The Hahn-Banach theorem is a topological result concerning separation of points from convex sets by hyperplanes, and the entire framework of Lagrange duality can be derived from this topological insight \citep{luenberger:book}. The Hahn-Banach theorem is also the basis for the universal representation theorem in deep neural networks \citep{DBLP:journals/mcss/Cybenko89}. Topological data analysis techniques, such as persistent homology, are playing an increasingly important role in different areas of machine learning \citep{DBLP:conf/sma/Edelsbrunner07, DBLP:journals/dcg/ZomorodianC05}. 

Graphical models have been extensively studied in artificial intelligence (AI), causal reasoning, machine learning (ML), physics, statistics, and many other fields \citep{koller:text,lauritzen:text, pearl:bnets-book,DBLP:journals/cacm/Pearl19}. Causal discovery \citep{spirtes:book} involves the construction of a causal model, for example a DAG graphical model structure and the specification of a probability model, from observational or experimental data. A broad family of models, ranging from Bayes networks \citep{pearl:bnets-book} on directed acyclic graphs (DAGs) to more recent variants, such as directed acyclic mixed graphs (ADMG) \citep{DBLP:conf/uai/Richardson09}, marginalized DAGs (mDAGs) \citep{mdag} and hyperedge-directed graphical models (HEDGs) \citep{hedge}, can be represented by finite Alexandroff spaces with different topological properties.  For example, a DAG model  imposes a $T_0$ topology on a finite Alexandroff space, which induces a partial ordering $\sigma$ on the variables $V$ in the model so that function $f_i$ determining the value of variable $X_{\sigma_i}$ in the model is measurable given the values of the previous variables $X_{\sigma_1}, \ldots, X_{\sigma_{i-1}}$. A directed acyclic mixed graph (ADMG) \citep{DBLP:conf/uai/Richardson09} and chain graphs \citep{DBLP:conf/uai/AnderssonMP96}, on the other hand, have both undirected and directed edges, which induce only a preordering on the set of variables. Marginalized DAGs (MDAGs) \citep{mdag} and HEDG \citep{hedge} models allow hyperedges between nodes, representing the effect of latent variables. These can be represented using topological constructions, such as non-Hausdorff cones $\mathbb{C}(X)$ or non-Hausdorff suspensions $\mathbb{S}(X)$ \citep{barmak}.

We propose a novel topological framework  for causal inference, building an initial topological representation of a partially ordered set (poset) from data, prior to building a probabilistic graphical model from the poset (see Figure~\ref{covid-diagram1}).  To that end, we represent posets using the algebraic topology of finite  Alexandroff spaces  \citep{alexandroff:1937,alexandroff:text,barmak,may}. Representing posets as topological spaces confers many computational advantages, such as the ability to 
combine multiple posets into a joint poset, and to use algebraic homotopy theory \citep{barmak,may} of finite topological spaces to significantly reduce the search space of possible structures. We show how a wide variety of graphical models, from chain graphs \citep{lauritzen:chain} to DAGs \citep{pearl:bnets-book}, can be topologically represented as posets in a finite Alexandroff space. We illustrate our approach using a real-world dataset of pancreatic cancer \citep{10.1093/biomet/asp023,ramon}.  Our primary goal is to illustrate how algebraic topology provides some powerful tools to design more scalable algorithms for structure discovery, which otherwise presents an intractable combinatorial search space. We also relate our approach to previous studies of causal discovery showing how to classify them using the concept of {\em intervention topology}  \citep{pearl:causalitybook,spirtes:book,DBLP:conf/uai/Eberhardt08,DBLP:journals/jmlr/HauserB12,DBLP:conf/nips/KocaogluSB17,MAOCHENG198415,prasad:aaai21}. \citep{pmlr-v108-bernstein20a} develop a greedy poset-based algorithm for learning DAG models, but do not exploit poset topological properties.

Our approach builds on a key technical innovation of using a topological representation of partially ordered sets, in between the original datasets and the final probabilistic or statistical graphical model. The principal reasons for explicitly modeling posets as topological spaces is that it allows us to exploit the rich algebraic theory of finite spaces \citep{barmak,may,mccord,stong} to our computational advantage. For example, in modeling cancer disease progression, in addition to obtaining crucial timing information about mutations from clinical datasets on tumors and their associated genotypes, we can also exploit the disease pathways \citep{Jones1801} that are known from medical research. Each source of information results in a different poset, which can be combined together using the algebraic topology theory of finite spaces. In addition, algebraic topology gives powerful tools for reducing the combinatorial search space of possible structures. We characterize homeomorphic equivalences among minimal poset models, and show homotopical equivalences sharply reduce the number of structures that need to be examined during structure discovery. In the case of pancreatic cancer, for example, we can reduce the search space of possible structures by three orders of magnitude.

The topology of posets  of graphical models based on the algebraic topology of finite  Alexandroff spaces  \citep{alexandroff:1937,alexandroff:text}.  We show how a wide variety of graphical models, from chain graphs \citep{lauritzen:chain} to DAGs \citep{pearl:bnets-book}, can be topologically embedded in a finite Alexandroff space. For a DAG ${\cal G}$, simply compute the unique transitive closure graph ${\cal G}_{tc}$, and define the open sets of the induced topological model $({\cal M}, \leq)$ by defining the open sets $U_x \subset {\cal T}$ as the ancestors of a node $x$ (including itself) in ${\cal G}_{tc}$. To revert from a topological model ${\cal M} = (X, {\cal T})$ to its DAG representation, form the Hasse diagram of the partial order defined by $x \leq y$ if $x \in U_y$. We  present a novel algorithmic paradigm for structure discovery as iterating between searching among topologically distinct structures and causally faithful structures. We show how this paradigm can be used to characterize many previous studies of causal discovery in terms of the concept of {\em intervention topology}, collections of subsets intervened on to determine directionality \citep{pearl:causalitybook,spirtes:book,DBLP:conf/uai/Eberhardt08,DBLP:journals/jmlr/HauserB12,DBLP:conf/nips/KocaogluSB17,MAOCHENG198415,prasad:aaai21}. This connection immediately suggests topological generalizations of these previous algorithms.

\section{Representing Causal DAGs as Finite Topological Spaces} 
\label{sec:top} 

We begin with brief review of basic point-set topology, and then give a succinct characterization of finite space topologies. Topology \citep{munkres:book} characterizes the abstract properties of arbitrary spaces that are equivalent under smooth deformations, usually represented as continuous invertible mappings called homeomorphisms. Formally, a general topological space $(X, \mathcal{U})$ is characterized by a base space $X$, along with a collection of ``open" sets ${\cal O}_i \subseteq \mathcal{U}$ closed under arbitrary union and {\em finite} intersection. Note the asymmetry in these restrictions. For example, if we define $X = \mathbb{R}$ to be the real line, and consider the open sets ${\cal O}_i = (-\frac{1}{i}, \frac{1}{i})$ to be the open intervals around $0$ for $i = 1, \ldots \infty$, then $\cap_{i=1}^\infty  (-\frac{1}{i}, \frac{1}{i}) = \{ 0 \}$, which is not an open set! 

Alexandroff \citep{alexandroff:1937,alexandroff:text} pioneered the study of the subclass of topological spaces that are closed under both arbitrary union and intersection. While our framework can be potentially be extended to the non-finite case, for simplicity, we will restrict our presentation in this paper to the case of finite Alexandroff spaces \citep{barmak,may,mccord,stong}. It is obvious to note that finite topological spaces $(X, \mathcal{U})$ are trivially closed under arbitrary unions and intersections, because there are only a finite number of open sets ${\cal O} \in \mathcal{U}$. However, what turns out to be surprising is that the particular construction used by Alexandroff in defining open sets applies even in the finite case, and results in spaces with surprising topological richness, even though they are finite.

\begin{definition} 
A finite Alexandroff topological space (or simply, finite space, in the remainder of the paper) $(X, \mathcal{U})$ is a finite set $X$ and a collection $\mathcal{U}$ of ``open" sets, namely subsets of $X$,  such that (i) $\emptyset$ (the empty set) and $X$ are in $\mathcal{U}$ (ii) Any union of sets in $\mathcal{U}$ is in $\mathcal{U}$ (iii) Any intersection of subsets in $\mathcal{U}$ is in $\mathcal{U}$ as well. 
\end{definition}

We will often refer to a finite space simply by its elements $X$, where the topology is left implicit, unless its character is important, when we will clarify it. The most common topologies on $X$ will be the {\em discrete} topology, where the collection of open sets $\mathcal{U}$ is just the powerset $2^{X}$, and the {\em trivial} topology $\mathcal{U} = \{ \emptyset, X \}$.

The major contribution of this paper is the use of a specific topological representation of partially ordered sets based on Alexandroff spaces \citep{alexandroff:1937,alexandroff:text}, who showed that finite topological spaces naturally defined preordered and partially ordered sets. Two classic papers by McCord \citep{mccord} and Stong \citep{stong} laid the foundations for much of the subsequent study of finite topological spaces. Detailed proofs of all the main theorems on finite topological spaces in this paper can be found in \citep{barmak,may,mccord,stong}.


Remarkably, a key idea that is implicit in many causal discovery  algorithms is the topological notion of {\em separability} \citep{cbn,DBLP:conf/nips/KocaogluSB17,acharya}, which intimately relates to the topology of the finite space. In order to construct a poset model from data, \citep{cbn} assume that the dataset $u: {\cal G} \rightarrow \mathbb{N}$ specifies the number of observations of each genotype $g$. For example,  Table~\ref{cancer-dataset} specifies the number of tumors that contain a specific set of gene mutations $g$. The {\em support} ${\cal S}(u)$ is the non-zero coordinates of $u$, namely the genotypes that occur in the data. A crucial assumption here is that the dataset $u$ {\em separates the events} $e$ and $f$ if there exists some genotype $g \in {\cal S}(u)$ such that $g \cap \{e, f \} \neq \emptyset$.  Viewed more abstractly, this notion of separability implies that the underlying space has the $T_0$ Kolmogorov topology. \citep{DBLP:conf/nips/KocaogluSB17} assume a separating set, which is essentially a restricted type of finite topological space. 

\begin{definition}
\label{t0}
The {\em neighborhood} of an element $x$ in a finite space $X$ is a subset $V \subset X$ such that $x \in U$ for some open set $U \subset V$. 
\begin{itemize} 
\item $X$ is a  Kolmogorov (or $T_0$) finite space $X$ if each pair of points $x, y \in X$ is distinguishable in the space, namely for each $x, y \in X$, there is an open set $U \in \mathcal{U}$ such that $x \in U$ and $y \notin U$. Alternatively, if $x \in U$  if and only if $y \in U,  \ \forall U \in \mathcal{U}$ implies that $x = y$. 
\item $X$ is a $T_1$ finite space if element $x \in X$ defines a closed set $\{ x \}$. 
\item $X$ is a $T_2$ finite space or a {\em Hausdorff} space if any two points have distinct neighborhoods.
\end{itemize}
\end{definition}

It turns out that $T_1$ finite spaces are not interesting since the only topology defined on them is the discrete (powerset) topology. The most interesting finite spaces are those equipped with the $T_0$ topology. 

\begin{lemma}\citep{may}
If $X$ is a $T_2$ space, then it is a $T_1$ space. If $X$ is a $T_1$ space, then it is a $T_0$ space.
\end{lemma}

The key concept that gives finite (Alexandroff) spaces its power is the definition of the minimal open basis. First, we introduce the concept of a basis in a topological space. 

\begin{definition}
A {\em basis} for the topological space $X$ is a collection $\mathcal{B}$ of subsets of $X$ such that 
\begin{itemize} 
\item For each $x \in X$, there is at least one $B \in \mathcal{B}$ such that $x \in B$.
\item If $x \in B' \cap B"$, where $B, B" \in \mathcal{B}$, then there is at least one $B \in \mathcal{B}$ such that $x \in B \subset B' \cap B"$. 
\end{itemize}
\end{definition}
The topology $\mathcal{U}$ {\em generated} by the basis $\mathcal{B}$ is the set of subsets $U$ such that for every $x \in U$, there is a $B \in \mathcal{B}$ such that $x \in B \subset U$. In other words, $U \in \mathcal{U}$ if and only if $U$ can be generated by taking unions of the sets in the basis $\mathcal{B}$. Now, we turn to giving the most important definition in Alexandroff spaces, namely the {\em unique minimal basis}. 

\begin{lemma}\citep{may}
Let $X$ be a finite Alexandroff space. For each $x \in X$, define the open set $U_x$ to be the intersection of all open sets that contain $x$. Define the relationship $\leq$ on $X$ by $x \leq y$ if $x \in U_y$, or equivalently, $U_x \subset U_y$ (where $x < y$ if the inclusion is strict). The open sets $U_x$ constitute a {\bf unique minimal  basis} $\mathcal{B}$ for $X$ in that if $\mathcal{C}$ is another basis for $X$, then $\mathcal{B} \subset \mathcal{C}$. Alternatively, define the closed sets $F_x = \{y \ | \ y \geq x \}$, which provide an equivalent characterization of finite Alexandroff spaces.\footnote{The minimal basic closed sets in a $T_0$ finite Alexandroff space correspond to the {\em ancestral sets} in a DAG graphical model.} 
\end{lemma}

Note that the relation $\leq$ defined above is a {\em preorder} because it is reflexive (clearly, $x \in U_x$) and transitive (if $x \in U_y$, and $y \in U_z$, then $x \in U_z$). However, in the special case where the finite space $X$ has a $T_0$ topology, then the relation $\leq$ becomes a partial ordering. This gives a topological way to model DAG models, which will play a crucial role in our framework. 

\begin{lemma}
A function $f: X \rightarrow Y$ from one finite space to another is continuous if and only if $f^{-1}(U)$ is open in $X$ if $U$ is open in $Y$.  
\end{lemma}

\begin{lemma}
If $X$ is a $T_2$ space, then it is a $T_1$ space. If $X$ is a $T_1$ space, then it is a $T_0$ space.
\end{lemma}

\begin{lemma}
A function $f: X \rightarrow Y$ between two finite spaces is continuous if and only if it is order-preserving, meaning if $x \leq x'$ for $x, x' \in X$, this implies $f(x) \leq f(x')$. 
\end{lemma}

 \begin{lemma}
 Let $x, y$ be two comparable points in a finite space $X$. Then, there exists a path from $x$ to $y$ in $X$, that is, a continuous map $\alpha: (0, 1) \rightarrow X$ such that $\alpha(0) = x$ and $\alpha(1) = y$. 
 \end{lemma}
 
 \begin{theorem}
 If $X$ is a finite topological space containing a point $y$ such that the only open (or closed) subset of $X$ containing $y$ is $X$ itself, then $X$ is contractible. In particular, the non-Hausdorff cone $\mathbb{C}(X)$ is contractible for any $X$. 
 \end{theorem}

{\bf Proof:} Let $Y = \{ * \}$ denote the space with a single element, $*$. Define the retraction mapping $r: X \rightarrow *$ by $r(x) = *$ for all $x \in X$, and define the {\em inclusion} mapping $i: Y \rightarrow X$ by $i(*) = y$. Clearly, $r \circ i  = id_{Y}.$ Define the homotopy $h: X \rightarrow I \rightarrow X$ by $h(x,t) = x$ if $t < 1$, and $h(x, 1) = y$. Then, $h$ is continuous, because for any open set $U$ in $X$, if $y \in U$, then clearly $U = X$ (as $X$ is the only open set containing $y$), and hence $h^{-1}(U) = X \times I$, which is open. If on the other hand, $y \notin U$, then $h^{-1}(U) = U \times [0, 1)$. It follows that $h$ is a homotopy $h \cong id_X = i \circ r$. \qed 

 \begin{lemma}
 If $X$ is an finite Alexandroff space, then $U_x$ is contractible. In particular, if $X$ has a unique maximal point or unique minimal point, then $X$ is contractible. 
 \end{lemma}

\subsection{Open sets induced by Causal DAGs} 

\begin{figure}[h!]
\centering
\begin{minipage}{1\textwidth}
\centering
\includegraphics[scale=0.3]{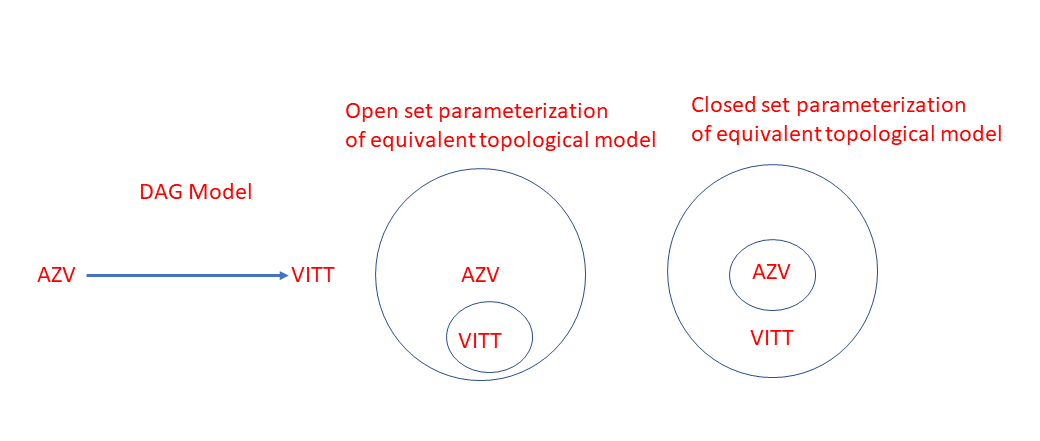}
\end{minipage}
\caption{A simple example of a causal relationship between {\bf AZV}, the adminstering of the AstraZeneca vaccine, and {\bf VITT}, the potentially vaccine-induced clotting of blood. On the left is a traditional DAG model of this causal relationship. In the middle, we give one characterization of this DAG as a finite space topological model, where the open sets of the model correspond to descendants of each node (including the node itself). On the right is an equivalent characterization in terms of closed sets, which correspond to ancestors of a node (including the node itself). }
\label{covid-diagram1}
\end{figure}

We use a motivational example of causal discovery for treatment of patients for the COVID pandemic using vaccines \citep{doi:10.1056/NEJMoa2104840,doi:10.1056/NEJMoa2104882}.
In the COVID vaccination causal discovery problem, we are given a finite set of $5$ variables $X = \{{\bf AZV}, {\bf PF4}, {\bf H}, {\bf BC}, {\bf VITT} \}$, defined as follows: 
\begin{itemize} 
\item {\bf AZV:} This variable represents the adminstration of the AstraZeneca vaccine. 

\item {\bf PF4:} A number of patients suffering from vaccine-induced abnormal blood clotting tested positive for heparin-induced platelet factor 4 (PF4). 

\item {\bf Gender:} Many of the patients who exhibited adverse effects to the Covid vaccine were disproportionately women, so gender may be a causal factor. 

{\bf HIT:} Heparin is a blood thinner used to prevent blood clots. Triggered by the immune system in response to heparin, HIT causes a low platelet count (thrombocytopenia).

\item {\bf VITT:} This variable denotes whether patients suffered from this rare vaccine-related variant of spontaneous heparin-induced thrombocytopenia that the authors of these studies referred to as vaccine-induced immune thrombotic thrombocytopenia. 
\end{itemize}
 Figure~\ref{covid-diagram1} illustrates a simple causal model for the COVID problem, represented both conventionally as a DAG, as well as two alternative representations of a finite topological space $(X, {\cal T})$, where $X$ is represented by the two variables shown, and ${\cal T}$ is a either a set of open sets, defined as the descendants of a node (including itself), or a set of closed sets comprised of the ancestors of a node (including itself).  In the open set parameterization, there is an arrow from node $x$ to node $y$ whenever $x \in {\cal O}_y$, that is, when the node $x$ is in the open set corresponding to node $y$. 

\begin{figure}[h!]
\centering
\begin{minipage}{1\textwidth}
\centering
\includegraphics[scale=0.30]{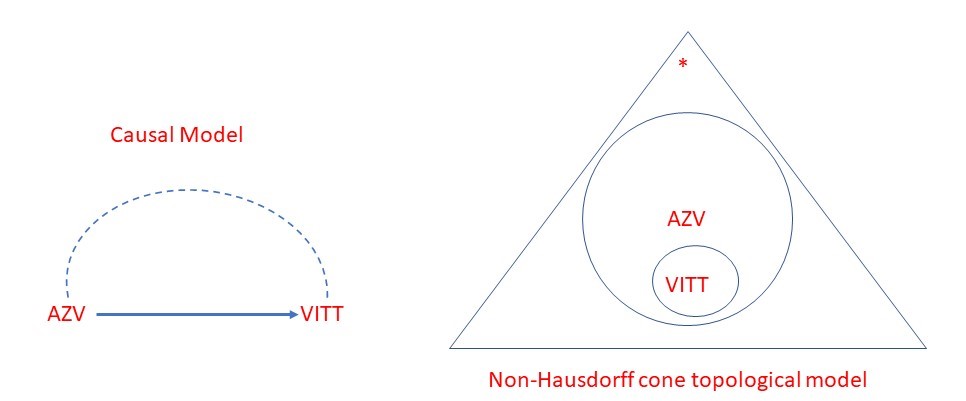}
\end{minipage}
\caption{A simple causal model with a latent variable correlating {\bf AZV}, the adminstering of the AstraZeneca vaccine, and {\bf VITT}, the potentially vaccine-induced clotting of blood. On the left is a traditional acyclic mixed directed graphical model (ADMG) of this causal relationship. On the right is the equivalent finite space topological model, where the open sets of the model correspond to descendants of each node (including the node itself), and the latent unobserved variable is represented by the element {\bf *} joined to the remainder of the model using a non-Hausdorff cone.}
\label{covid-diagram2}
\end{figure}

 Figure~\ref{covid-diagram2} shows how latent variables are typically modeled in causal reasoning using acyclic mixed directed graphs (ADMGs), where the latent variable is represented by a dashed undirected edge connecting the two observable variables. Such latent variables can be captured in our finite space topological framework by the use of the non-Hausdorff cone construction, which is one of several ways of connecting two topological spaces. Recall the non-Hausdorff cone merging of topological space $X$ with $\{ * \}$ yields the new space $\mathbb{C}(X)$, whose open sets are now ${\cal O}_{\mathbb{C}(X)} = {\cal O}_X \cup \{X \cup \{ * \} \}$. 

To take a real-world example, recently two studies were published in the New England Journal of Medicine that described patients in Austria, Germany and Norway who developed an unexpected blood clotting disorder in reaction to their first dose of the AstraZeneca/Oxford COVID-19 vaccine \citep{doi:10.1056/NEJMoa2104840,doi:10.1056/NEJMoa2104882}. Understanding causal pathways in such problems requires modeling the effects of tens of thousands of discrete and continuous variables, from the administration of the vaccine, heparin-induced platelet factors like PF4, thromocytopenia (blood clots) and its various causes, and the entire previous medical history of the patient. Clinicians have to juggle through all these factors in describing a potential treatment (e.g., should heparin be given to a patient?).

\begin{figure}[h!]
\centering
\begin{minipage}{1\textwidth}
\centering
\includegraphics[scale=0.30]{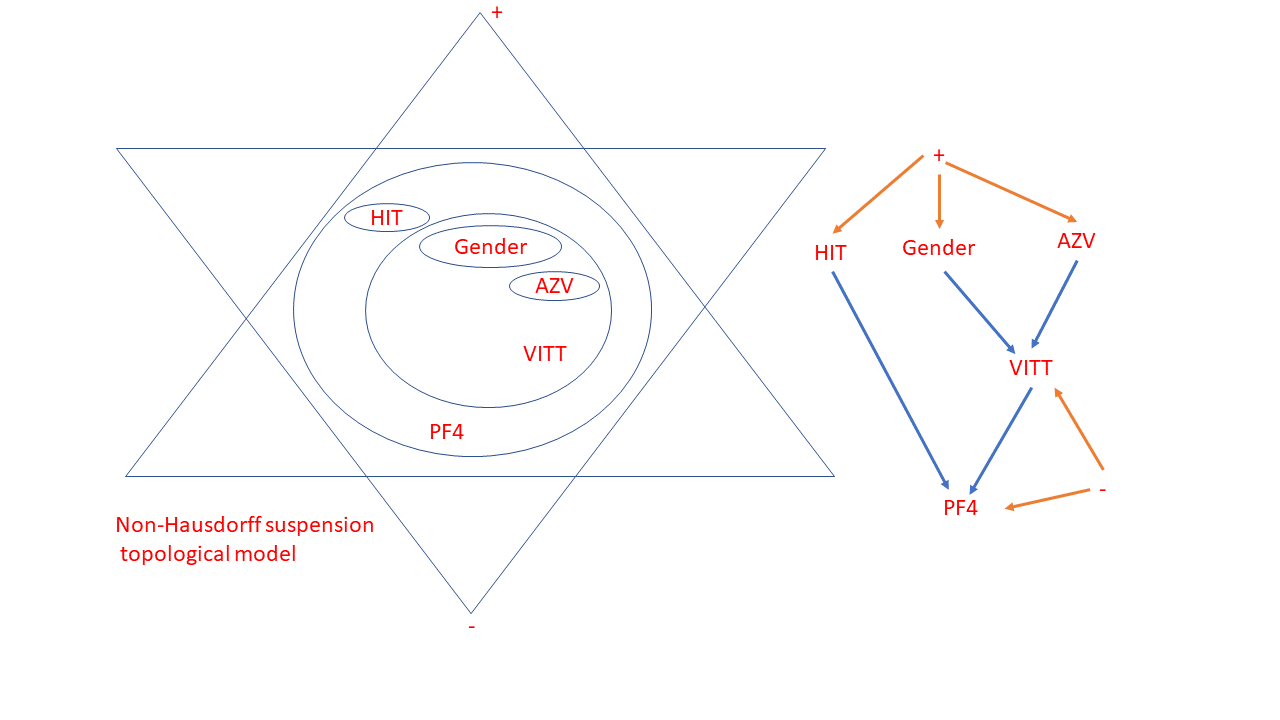}
\end{minipage}
\caption{A directed graph with hyper-edges (HEDG) represented as a finite topological space. Open sets associated with variables are shown as circles or ellipsoids. The ${\bf +}$ and ${\bf -}$ nodes are two latent factors connected by hyper-edges to observables. The {\em non-Hausdorff suspension} $\mathbb{S}$ defined in Section~\ref{sec:top} permits composing the hyper-edge latent variable model with the topological model induced by the DAG over observables. }
\label{HEDG-model}
\end{figure}

\subsection{Connectivity in Topological Spaces}
\label{connected} 

As mentioned above, every concept in a topological space must be defined in terms of the open (or closed) set topology, and that includes (path) connectivity. The crucial idea here is that connectivity is defined in terms of a continuous mapping from the unit interval $I = (0,1)$ to a topological space $X$. Remarkably, the upshot of this construction is that every graph-theoretic concept in causal models, e.g. separation and conditional independence, can be translated into properties of open sets in the topology. 

\begin{lemma}
A function $f: X \rightarrow Y$ between two finite spaces is continuous if and only if it is order-preserving, meaning if $x \leq x'$ for $x, x' \in X$, this implies $f(x) \leq f(x')$. 
\end{lemma}

\begin{definition}
We call two points $x, y \in X$ {\bf comparable} if there is a sequence of elements $x_0, \ldots, x_n$, where $x_0 = x, x_n = y$ and for each pair $x_i, x_{i+1}$ either $x_i \leq x_{i+1}$ or $x_i \geq x_{i+1}$. A {\bf fence} in $X$ is a sequence $x_0, x_1, \ldots, x_n$ of elements such that any two consecutive elements are comparable. $X$ is {\bf order connected} if for any two elements $x, y \in X$, there exists a fence starting in $x$ and ending in $y$. 
\end{definition}
 
 \begin{lemma}
 Let $x, y$ be two comparable points in a finite space $X$. Then, there exists a path from $x$ to $y$ in $X$, that is, a continuous map $\alpha: (0, 1) \rightarrow X$ such that $\alpha(0) = x$ and $\alpha(1) = y$. 
 \end{lemma}
 
 \begin{lemma}
 Let $X$ be a finite space. The following are equivalent: (i) $X$ is a connected topological space. (ii)  $X$ is an order-connected topological space (iii) $X$ is a path-connected topological space. 
 \end{lemma}
 
 To illustrate the notion of connectivity, Table~\ref{enum-t0} gives examples of connected and disconnected finite space topologies for a small three element space. 
 
 \subsection{Separation and Conditional Independence in Topological Spaces}
\label{ci-top} 

We now give a purely topological characterization of separation and conditional independence in finite topological spaces, which draw upon equivalent notions in graphical models \citep{lauritzen:chain,pearl:bnets-book}, but are defined with respect to the open sets of the topology. 

\begin{definition}
Given a connected finite space $X$, with an induced (pre,partial) ordering $\leq$, and subsets $U, V, Z \subset X$, the subset $U$ is {\bf topologically blocked or d-separated} from $V$ given $Z$, if for every fence from an element $u \in V$ to an element $v \in V$, the following conditions hold: 
\begin{enumerate}
    \item The fence $x_0, x_1, \ldots, x_n$, where $x_0 = u$ and $x_n = v$ is such that every consecutive pair of elements is of the form $x_i \leq x_{i+1}$ or $x_{i+1} \leq x_i$, and some element $x_k \in Z$ (this condition is equivalent to stating that all the edges are of the form $x_i \rightarrow x_{i+1}$ or $x_{i} \leftarrow x_{i+1}$). 
    \item The fence $x_0, x_1, \ldots, x_n$, where $x_0 = u$ and $x_n = v$ is such that for every {\bf collider} in the fence, namely a triple of elements $x_i, x_{i+1}, x_{i+2}$ is such that $x_i \leq x_{i+1}$, and $x_{i+2} \leq x_{i+1}$, it holds that $U_{x_{i+1}} \cap Z = \emptyset$. This condition is the topological restatement of the standard collider condition in graphical models, where for a path to be blocked, no collider or any of its descendants can be in the conditioning set.
\end{enumerate}
\end{definition}

\begin{definition}
Given a finite space $X$, and subsets $U, V, Z \subset X$, $U$ is {\bf topologically conditionally independent} (TCI) of $V$ given $Z$ if and only if every fence from an element $u \in U$ to an element $v \in V$ is topologically blocked with respect to the conditioning set $Z$. 
\end{definition}

\section{Stable and Solvable Causal Models over Finite Topological Spaces} 
\label{wif}

We now precisely define {\em stable} and {\em solvable} causal models, including both their Alexandroff topological structure, and a decomposable product probability measure constituting the parameters of the model.  We impose the condition that the product probability measure respect the underlying Alexandroff topology, namely its open (or closed) sets, and requirement of a particular factorization is translated into a requirement of a particular topology.  Our formulation is related to the {\em intrinsic} model of decision making \citep{witsenhausen:1975}, which was recently adapted to causal inference \citep{cif}. Neither of these investigated Alexandroff topologies.

\begin{definition}
A {\bf finite space causal model} is defined as ${\cal M} = (U_\alpha, {\cal F}_\alpha, {\cal I}_\alpha, (\Omega, {\cal B}, P))$, where $\alpha \in X$, a finite Alexandroff space topology. $U_\alpha$ is a non-empty set that defines the range of values that variable $\alpha$ can take. ${\cal F}_\alpha$ is a sigma field (or algebra) of measurable sets for variable $\alpha$. The triple $(\Omega, {\cal B} , P)$ is a probability space, where ${\cal B}$ is a sigma field of measurable subsets of sample space $\Omega$. The {\bf information field} ${\cal I}_\alpha \subset {\cal F}$ represents the ``receptive field" of an element $\alpha \in X$, namely the set of other elements $\beta \in X$ whose values $\alpha$ must consult in determining its own value. We impose the restriction that the information field ${\cal I}_\alpha$ respect the Alexandroff topology on $X$, so that ${\cal I}_\alpha \subset {\cal F}(U_\alpha)$, where $U_\alpha$ is the minimal basic open set associated with element $\alpha \in X$.
\end{definition}

Following structural causal models \citep{pearl:bnets-book}, we can decompose the elements of the topological space into disjoint subsets $X = U \sqcup V$, where $U$ represents ``exogenous" variables that have no parents, namely $\alpha$ is exogenous precisely when ${\cal I}_\alpha \subset {\cal F}(\emptyset)$, and $V$ are ``endogenous" variables whose values are defined by measurable functions over exogenous and endogenous variables. Note that the probability space can be defined over the ``exogenous" variables $\alpha \in U$, in which case it is convenient to attach a local probability space $(\Omega_\alpha, {\cal B}_\alpha, P)$ to each exogenous variable, where ${\cal B}_\alpha \subset {\cal B}$. We define conditional independence with respect to the induced information fields over the open sets of the Alexandroff space. 

\begin{definition}
Given the induced probability space over information fields in a topological finite space, a {\em stochastic basis} is a sequence of information fields ${\cal G} = {\cal I}_1, \ldots, {\cal I}_n$ such that for $1 \leq i \leq n-1, {\cal I}_i \subset {\cal I}_{i+1}$, and $\cup_{1=1}^n {\cal I}_i = {\cal F}$. Two such sequences ${\cal G}_1$ and ${\cal G}_2$ are {\bf conditionally independent} given the base sigma algebra ${\cal F}$, if for all subsets $A \in {\cal G}_1$, $B \in {\cal G}_2$, it follows that
$P(A \ B | {\cal F} )= P(A | {\cal F}) P(B | {\cal F})$. 
\end{definition}

\begin{definition}
The {\bf decision field} $U = \prod_{\alpha \in X} U_\alpha$ defines the space of all possible values of the variables in the finite space causal model, where the cartesian product is interpreted as a map $u: X \rightarrow \cup_{\alpha \in X} U_\alpha$ such that $u(\alpha) \equiv u_\alpha \in U_\alpha$.  
\end{definition}

\begin{definition}
For any subset of elements $B \in X$, let $P_B$ denote the {\em projection} of the product $\prod_\alpha U_\alpha$ upon the product $\prod_{\beta \in B} U_\beta$, that is $P_B(u)$ is simply the restriction of $u$ to the domain $B$. 
\end{definition}

\begin{definition}
The {\bf product sigma field} is defined as $\prod_{\alpha \in B} {\cal F}_\alpha$ over $\prod_{\alpha \in B} U_\alpha$, where ${\cal F}(B)$ is the smallest sigma-field such that $P_B$ is measurable. Note that if $B_1 \subset B_2$, then ${\cal F}_{B_1} \subset {\cal F}_{B_2}$. The finest sigma-field ${\cal F}(X) = \prod_{\alpha \in X} {\cal F}_\alpha$.  
\end{definition}

\begin{definition}
A finite space causal model ${\cal M}$ is {\bf causally faithful} with respect to the probability distribution $P$ over ${\cal M}$ if every conditional independence in the topology, as defined in Section~\ref{ci-top}, is satisfied by the distribution $P$, and vice-versa, every conditional independence property of the $P$ is satisfied by the topology. 
\end{definition}

We can now formally define what it means to ``solve" a causal finite space model ${\cal M}$. We impose the requirement that each variable $\alpha \in X$ must compute its value using a function measurable on its own information field. 

\begin{definition}
Let the policy function $f_\alpha$ of each element $\alpha \in X$ be constrained so that $f_\alpha: U \times \Omega \rightarrow U_\alpha$ is measurable on the product sigma field ${\cal I}_\alpha \times {\cal B}_\alpha$, namely $f_\alpha^{-1}({\cal F}_\alpha) \subset {\cal I}_\alpha \times {\cal B_\alpha}$. 
\end{definition}

\begin{definition}
\label{property-sm1}
The finite space causal model ${\cal M} = (X, U_\alpha, {\cal F}_\alpha, {\cal I}_\alpha)$ is {\bf measurably solvable} if for every $\omega \in \Omega$, the closed loop equations $P_\alpha(u) = f_\alpha(u, \omega)$ have a unique solution for all $\alpha \in X$, where for a fixed $\omega \in \Omega$, the induced map ${\cal M}^\gamma: \Omega \rightarrow U$ is a measurable function from the measurable space $(\Omega, {\cal B})$ into $(U, {\cal F})$. 
\end{definition}

\begin{definition}
\label{property-sm2}
The finite space causal model ${\cal M} = (X, U_\alpha, {\cal F}_\alpha, {\cal I}_\alpha)$ is {\bf stable} if for every $\omega \in \Omega$, the closed loop equations $P_\alpha(u) = f_\alpha(u, \omega)$ are solvable by a fixed constant ordering $\Xi$ that does not depend on $\omega \in \Omega$. 
\end{definition}

Measurably solvable models capture the corresponding property in a structural causal model $(U,V,F,P)$, which states that for any fixed probability distribution $P$ defined over the exogenous variables $U$, each function $f_i$ computes the value of variable $x_i \in V$, given the value of its parents $Pa(x_i)$ uniquely as a function of $u \in U$. This allows defining the induced distribution $P_u(V)$ over exogenous variables in a unique functional manner depending on some particular instantiation of the random exogenous variables $U$. Stable models are those where the ordering of variables is fixed.  We now extend the notion of {\em recursive} causal models in DAGs \citep{pearl:causalitybook} to finite topological spaces. 

\begin{definition}
The finite space causal model ${\cal M} = (X, U_\alpha, {\cal F}_\alpha, {\cal I}_\alpha)$ is a {\bf recursively causal} model if there exists an ordering function $\psi: X \rightarrow \Xi_n$, where $\Xi_n$ is the set of all injective (1-1) mappings of $(1, \ldots, n)$ to the set $X$, such that for any $ 1 \leq k \leq n$, the information field of variable $\alpha_k$ in the ordering $\Xi_n$ is contained in the joint information fields of the variables preceding it:
\begin{equation}
    {\cal I}_{\alpha_k} \subset {\cal F}(\alpha_1, \ldots, \alpha_{k-1})
\end{equation}
\end{definition}


\begin{definition}
A {\bf causal intervention}  do$(\beta$=$u_\beta)$ in a finite space topological model ${\cal M} = (X, U_\alpha, {\cal F}_\alpha, {\cal I}_{\alpha})$ is defined as the submodel ${\cal M}_\beta$ whose information fields ${\cal I}_\alpha$ are exactly the same as in $M$ for all elements $\alpha \neq \beta$, and the information field of the intervened element $\beta$ is defined to be ${\cal I}_\beta \subset {\cal F}(\emptyset) \times {\cal B}_\beta$. Note that since the only measurable function on ${\cal F}(\emptyset)$ is the constant function, whose value depends on a random sample space element $\omega \in \Omega_\beta$, this generalizes the notion of causal intervention in DAGs, where an intervened node has all its incoming edges deleted. \footnote{Our definition of causal intervention differs from that proposed in causal information fields \citep{cif}, where additional intervention nodes were added to the model.} 
\end{definition}

\subsection{Embedding Causal Graphical Models into Finite Topological Spaces} 
\label{gm}

We now explain how to construct faithful topological embeddings of causal graphical models. The following lemma plays a fundamental role in constructing topologically faithful embeddings of graphical models. 

\begin{lemma} \citep{barmak,may}
\label{preorder}
A preorder $\{X, \leq \}$ determines a topology $\mathcal{U}$ on space $X$ with the basis given by the collection of open sets $U_x = \{y \ | \  y \leq x \}$. It is referred to as the {\bf order topology} on $X$. The space $(X, \mathcal{U})$ is a $T_0$ space if and only if $(X, \leq)$ is a partially ordered set (poset). As before, we can alternatively characterize finite space topologies by the closed sets $F_x = \{y \ | \ y \geq x \}$. 
\end{lemma}
The unique minimal basis gives us a way of characterizing whether or not a finite space has $T_0$ topology. 

\begin{lemma} \citep{may}
Two elements $x, y \in X$ have the same neighborhoods if and only if $U_x = U_y$. Thus, a finite space $X$ has $T_0$ topology if and only if $U_x = U_y$. 
\end{lemma} 

{\bf Proof:} If $x$ and $y$ have the same neighborhoods, then trivially $U_x = U_y$. Conversely, if $U_x = U_y$, then if $x \in U$ for some open set $U$, then $U_y = U_x \subset U$ (recall that $U_x$ is the intersection of all sets that contain $x$), and hence $y \in U$. Similarly, if $y \in U$, the same argument shows $x \in U$. Thus, $x$ and $y$ have the same neighborhoods. $\qed$

We state three theorems that show how to reduce several popular causal models into their faithful topological embedding. The same construction can be followed for all the other models in the literature as well. We focus on embedding the topology of a graphical model, leaving aside the parametric specification of a probability measure on the model (which we discuss in more depth in the appendix). 

\begin{theorem}
Every causal DAG  graphical structure $G = (V,E)$ defines a finite $T_0$ Alexandroff topological space with a partial ordering.  
\end{theorem}

{\bf Proof:} Define the elements of the topology $X = V$, the vertices representing the variables of the DAG $G = (V,E)$. Construct the transitive closure $G_{tc}$ of the DAG $G$.  Define the partial ordering $x \leq y$ in the topological space if the variable $x$ is a descendant of $y$ in $G_{tc}$.  Define the open sets of $X$ as $U_y = \{x | x \leq y \}$.  $\qed$

\begin{theorem}
Every  chain  graph \citep{lauritzen:chain} structure $G = (V,E)$ defines a finite Alexandroff topological space $X$ with a $\leq$ preordering. 
\end{theorem}

{\bf Proof:} Once again, define $X$ as the variables in the chain graph. Recall that in a chain graph $G = (V, E)$, two nodes $x$ and $y$ are connected by an edge that is either directed, so $x \leftarrow y$ or $x \rightarrow y$, or there is an undirected edge  $x - y$ between them. Define the ordering $x \leq y$ on the topology $X$ if and only if there exists a path from $x$ to $y$ such that every comparable pair of nodes on this path is either of the form $x_k \leftarrow x_{k+1}$ or alternatively $x_{k} - x_{k+1}$. This ordering $\leq$ on $X$ is a preordering, and hence defines a general Alexandroff finite space topology. Define the open sets of $X$ as $U_x = \{y | y \leq x \}$.  $\qed$

\begin{theorem}
Every  mDAG graphical model \citep{mdag} or HEDG hyper-edge directed graphical model \citep{hedge} $G = (V, E, H)$, where $H$ is a set of hyper-edges, represented by an abstract simplicial complex, can be represented by a finite Alexandroff topological space $X$ with a $\leq$ preordering. 
\end{theorem}

{\bf Proof:} Define the space $X$ by the variables in the graphical model $G = (V, E, H)$. For the observable edges represented by $E$, we follow the same construction as in DAG models described above. Note the hyper-edges $h \in H$ in effect represent an abstract simplicial complex. For example, in Table~\ref{enum-t0}, the $D_3$ discrete topology on $X = \{a, b, c\}$ can be represented as an mDAG model where the three observable variables $a, b, c$ are connected only through one latent variable, whose effect on the observable variables is manifested by the hyper-edge that constitutes a simplicial complex ${\cal C }$ defined by the non-empty power set of $X$. This simplicial complex $ {\cal C}$ can be modeled as a non-Hausdorff cone $\mathbb{C}$ between the latent variable and the open set topology of the observable variables (see Section~\ref{prod-top} below). $\qed$

\begin{table}
 \caption{Examples of finite Alexandroff spaces on $X =\{a, b, c \}$, and their equivalent graphical models. A proper open set is any set other than $\emptyset$ or $X$ (which are in any topology). $P_n$ is a topology on a set of size $n$ with only one proper open set. $D_n$ is the discrete topology over $n$ elements. $P_{m,n}$ are topologies where the proper open sets are all non-empty subsets of a subset of size $m$. The $\cong$ equivalence relation  is homotopy equivalence. See text for explanation.}
 \centering
 \begin{small}
  \begin{tabular}{|c|c|c|c|c|} \hline 
Proper Open Sets & Name & $T_0$? & Connected? & Equivalent graphical model \\ \hline 
All & $D_3$ & yes & no & HEDG (hyper-edge over (a,b,c)) \citep{hedge} \\ \hline 
$b, c$ & & yes & yes & DAG $b \rightarrow a$, $c \rightarrow a$ (collider over a) \\ \hline
$a, b, (a,b)$  & $P_{2,3} \cong \mathbb{C}D_2$  & yes &  yes & Chain graph: $a \rightarrow c, b \rightarrow c, a - b$ \\ \hline 
$a, b, (a,b), (b,c)$ & $D_1 \bigsqcup P_2$  & yes &  no & DAG with node $a$ disconnected, $b \rightarrow c$\\ \hline 

$a$  & $P_3$  & no &  yes &  Chain graph: $a \rightarrow b$, $a \rightarrow c$, and $b - c$\\ \hline
\end{tabular}
\end{small}
\label{example-enum}
\end{table}

\subsection{Combining Poset Models}

 A crucial strength of our topological framework is the ability to combine two topological spaces $X$ and $Y$ into a new space, which can generate a rich panoply of models. Here are a few of the myriad ways in which topological spaces can be combined \citep{munkres:book}. Table~\ref{example-enum} illustrates some of these ways of combining spaces for a small finite space $X$ comprised of just three elements. 

\begin{itemize} 
\label{prod-top}
\item Subspaces: The {\bf subspace} topology on $A \subset X$ is defined by the set of all intersections $A \cap U$ for open sets $U$ over $X$. 
\item Quotient topology: The {\bf quotient topology} on $U$ defined by a surjective mapping $q: X \rightarrow Y$ is the set of subsets $U$ such that $q^{-1}(U)$ is open on $X$. 
\item Union: The {\bf topology of the union} of two spaces $X$ and $Y$ is given by their disjoint union $X \bigsqcup Y$, which has as its open sets the unions of the open sets of $X$ and that of $Y$. 
\item Product of two spaces: The {\bf product topology} on the cartesian product $X \times Y$ is the topology with basis the ``rectangles" $U \times V$ of an open set $U$ in $X$ with an open set $V$ in $Y$.
\item Wedge sum of two spaces: The {\em wedge sum} is the ``one point" union of two ``pointed" spaces $(X, x_o)$ with $(Y, y_o)$, defined by $X \bigvee Y / x_0 \sim y_0$, the quotient space of the disjoint union of $X$ and $Y$, where $x_0$ and $y_0$ are identified. 
\item Smash product: The {\em smash product} topology is defined as the quotient topology $X \bigwedge Y = X \times Y / X \bigvee Y$. 
\item Non-Hausdorff cone: The {\bf non-Hausdorff cone} of topological space $X$ with $Y = \{ * \}$ yields the new space $\mathbb{C}(X)$, whose open sets are now ${\cal O}_{\mathbb{C}(X)} = {\cal O}_X \cup \{X \cup \{ * \} \}$. 
\item Non-Hausdorff suspension: The {\bf  non-Hausdorff suspension} of topological space $X$ with $Y = \{+, - \}$ yields the new space $\mathbb{S}(X)$, whose open sets are now ${\cal O}_{\mathbb{C}(X)} = {\cal O}_X \cup \{X \cup \ \{+, - \} \}$. 
\end{itemize} 

\subsection{Homeomorphisms and Homotopical Equivalences} 
\label{homotopy} 

\begin{table}
 \caption{Enumeration of Alexandroff finite space topologies (see \citep{may}).}
 \centering
  \begin{tabular}{|l|l|l|l|l|} \hline 
n & Distinct & Distinct $T_0$ & Inequivalent & Inequivalent $T_0$ \\ \hline 
1 & 1 & 1 & 1 & 1 \\ \hline 
2 & 4  & 3 &  3 &  2 \\ \hline 
3  & 29  & 19 &  9  & 5 \\ \hline 
4  & 355  & 219 &  33  & 16\\ \hline
5  & 6942  & 4231 &  139 &  63\\ \hline
6  & 209,527  & 130,023  & 718  & 318\\ \hline
7  & 9,535,241  & 6,129,859  & 4,535  & 2,045\\ \hline
8  & 642,779,354  & 431,723,379  & 35,979  & 16,999\\ \hline
9  & 63,260,289,423  & 44,511,042,511  & 363,083  & 183,231\\ \hline
10  & 8,977,053,873,043  & 6,611,065,248,783  & 4,717,687  & 2,567,284\\ \hline
\end{tabular}
\label{top-enum-t0}
\end{table}

Like the number of possible DAG structures, the number of possible finite space topologies grows extremely rapidly. However, there are powerful tools in algebraic topology, such as homotopies, which characterize equivalences among spaces \citep{munkres:book}. In particular, Table~\ref{top-enum-t0} shows that exploiting homotopical inequivalences, we can save over three orders of magnitude in searching for an appropriate poset model over naive search. Given that evolutionary processes, such as pancreatic cancer, may potentially involve multiple thousands of elements (genes), the savings may be very significant. Of course, it is crucial to combine the savings from domain knowledge, as provided in Table~\ref{cancer-dataset} and Table~\ref{pathways} with that provided by efficient enumeration of poset topologies under homemorphic equivalences. 

\begin{definition}
 A topological space $X$ is {\bf contractible} if the identity map $id_X: X \rightarrow X$ is homotopically equivalent to the constant map $f(x) = c$ for some $c \in X$. 
 \end{definition}
 
For example, any convex subset $A \subset \mathbb{R}^n$ is contractible. Let $f(x) = c, c \in A$ be the constant map. Define the homotopy $H: A \times I \rightarrow X$ as equal to $H(x,t) = t c + (1 - t) x$. Note that at $t=0$, we have $H(x,0) = x$, and that at $t=1$, we have $H(x,1) = c$, and since $A$ is a convex subset, the convex combination $t c + (1 - t) x \in A$ for any $t \in [0,1]$. 
 
 \begin{theorem}
 If $X$ is a finite topological space containing a point $y$ such that the only open (or closed) subset of $X$ containing $y$ is $X$ itself, then $X$ is contractible. In particular, the non-Hausdorff cone $\mathbb{C}(X)$ is contractible for any $X$. 
 \end{theorem}

{\bf Proof:} Let $Y = \{ * \}$ denote the space with a single element, $*$. Define the retraction mapping $r: X \rightarrow *$ by $r(x) = *$ for all $x \in X$, and define the {\em inclusion} mapping $i: Y \rightarrow X$ by $i(*) = y$. Clearly, $r \circ i  = id_{Y}.$ Define the homotopy $h: X \rightarrow I \rightarrow X$ by $h(x,t) = x$ if $t < 1$, and $h(x, 1) = y$. Then, $h$ is continuous, because for any open set $U$ in $X$, if $y \in U$, then clearly $U = X$ (as $X$ is the only open set containing $y$), and hence $h^{-1}(U) = X \times I$, which is open. If on the other hand, $y \notin U$, then $h^{-1}(U) = U \times [0, 1)$. It follows that $h$ is a homotopy $h \cong id_X = i \circ r$. \qed 

The following lemma is of crucial importance in Section~\ref{sec:algm}, where we will define {\em beat points}, elements of a topological space that can be removed, reducing model size. 

\begin{definition}
A point $x$ in a finite Alexandroff topological space $X$ is {\bf maximal} if there is no $y > x$, and {\bf minimal} if there is no $y < x$. 
\end{definition}
 
 \begin{lemma}
 If $X$ is an finite Alexandroff space, then $U_x$ is contractible. In particular, if $X$ has a unique maximal point or unique minimal point, then $X$ is contractible. 
 \end{lemma}

\begin{definition}
Let $f, g: X \rightarrow Y$ be two continuous maps between finite space topologies $X$ and $Y$. We say $f$ is {\em homotopic} to $g$, denoted as $f \cong g$ if there exists a continuous map $h: X \times [0,1] \rightarrow Y$ such that $h(x, 0) = f(x)$ and $h(x, 1) = g(x)$. In other words, there is a smooth ``deformation" between $f$ and $g$, so we can visualize $f$ being slowly warped into $g$. Note that $\cong$ is an equivalence relation, since $f \cong f$ (reflexivity), and if $f \cong g$, then $g \cong f$ (symmetry), and finally $f \cong g, g \cong h \ \ \implies f \cong h$ (transitivity). 
\end{definition}
 
 \begin{definition}
 A map $f: X \rightarrow Y$ is a {\em homotopy equivalence} if there exists another map $g: Y \rightarrow X$ such that $g \circ f \cong id_X$ and $f \circ g \cong id_Y$, where $id_X$ and $id_Y$ are the identity mappings on $X$ and $Y$, respectively. 
 \end{definition}

\section{Algorithms for Learning Causal Posets} 
\label{sec:algm} 

In this section, we describe a number of algorithms for constructing causal poset models. Our approach is intended to highlight the important role played by topological constraints, which were implicit in many previous studies. We use the domain of cancer genomics to illustrate how topological constraints on datasets makes it possible to efficiently learn the structure and parameters of a causal model from observational data \citep{10.1093/biomet/asp023,cbn,BEERENWINKEL2006409,gerstung}.  A greedy algorithm for learning causal poset models is described in \citep{pmlr-v108-bernstein20a}, but it does not fully exploit the algebraic topology of posets for efficient enumeration. We show that the space of causal structures can be significantly pruned by exploiting the algebraic topology of finite spaces, in particular using homeomorphisms among topologically equivalent finite space models. 

\subsection{Intervention Topologies}

First, we show that many previous studies of causal discovery from interventions, including the {\em conservative} family of intervention targets \citep{DBLP:journals/jmlr/HauserB12}, path queries \citep{DBLP:conf/nips/BelloH18a}, and {\em separating systems} of finite sets or graphs.  \citep{DBLP:conf/uai/Eberhardt08,DBLP:journals/jmlr/HauserB12,DBLP:conf/nips/KocaogluSB17,sepsets,MAOCHENG198415} can all be viewed as imposing an intervention topology. Table~\ref{intervention-top} classifies a few previous studies in terms of the induced intervention topology.

\begin{table}[b]
 \caption{Intervention topologies of some previous causal discovery methods.}
 \centering
 \begin{small}
  \begin{tabular}{|c|c|c|c|c|} \hline 
Intervention Sets & Topology & Intervention Class & Causal Structure & Reference \\ \hline 
$\{\emptyset, X \}$  & Trivial & Observations & Conjunctive Bayes Network & \citep{BEERENWINKEL2006409} \\ \hline 
$I \sqcup X \setminus I$  & Disconnected & All subsets &  DAG & \citep{DBLP:conf/uai/Eberhardt08} \\ \hline 
$F_x = \{ x \}$ & $T_1$ & Single-node & Tree graphs  & \citep{shpitser} \\ \hline 
 Strong separating sets & $T_2$ & Restricted open sets & Latent DAG & \citep{DBLP:conf/nips/KocaogluSB17} \\ \hline
$\{ x \}$ &  $T_1$ & Path queries & Transitive DAG & \citep{DBLP:conf/nips/BelloH18a} \\ \hline 
 Minimal elements & Restricted open sets & Leaf queries & Tree Graphs & \citep{prasad:aaai21} \\ \hline
\end{tabular}
\end{small}
\label{intervention-top}
\end{table}
 
 If no experiments are allowed, the intervention topology is simply the trivial topology. If $X = I \sqcup X \setminus I$, where $I$ is the intervention target \citep{DBLP:conf/uai/Eberhardt08}, then the intervention topology is disconnected. \citep{shpitser} study single node interventions, which can be viewed as a $T_1$ intervention topology where singleton sets are closed. \citep{DBLP:journals/jmlr/HauserB12} introduce the idea of {\em conservative} family ${\cal I}$ of intervention targets, meaning a family of (open) subsets of variables in a causal model such that for every variable $x \in X$, there exists an $I \subset {\cal I}$ such that $x \notin I$. This is closely related to the idea of Alexandroff topologies where elements have distinguishable neighborhoods, and each intervention target $I \in {\cal I}$ defines a neighborhood. A very related notion is that of {\em separating systems} of finite sets as intervention targets  \citep{DBLP:conf/uai/Eberhardt08,DBLP:journals/jmlr/HauserB12,DBLP:conf/nips/KocaogluSB17}. or separating systems of graphs \citep{MAOCHENG198415,DBLP:journals/corr/hauser-arxiv}. \citep{DBLP:conf/nips/KocaogluSB17} used {\em antichains}, a partitioning of a poset into subsets of non-comparable elements. \citep{DBLP:conf/nips/BelloH18a} use path queries, which can be viewed as chains. Finally, \citep{prasad:aaai21} used leaf queries on tree structures, where none of the interior nodes can be intervened on. 

We first introduce the notion of a {\em separating system}, which is a special case of the $T_0$ topology separation axiom of finite Alexandroff spaces. 

\begin{definition}
A {\em separating system} on a finite set $X$ is a collection of subsets $\{\mathcal{U}_1, \ldots, \mathcal{U}_m \}$ such that for every pair of elements $x, y \in X$, there is a set $\mathcal{U}_k$ such that either $x \in \mathcal{U}_k, y \notin \mathcal{U}_k$ or alternatively, $x \notin \mathcal{U}_k, y \in \mathcal{U}_k$. An $(m,n)$ {\em strongly separating system} is a pair of sets $\mathcal{U}_i, \mathcal{U}_j$ such that $x \in \mathcal{U}_i, y \notin \mathcal{U}_j$ and $x \notin \mathcal{U}_i, y \in \mathcal{U}_j$. 
\end{definition}

\begin{definition}
Given a finite space Alexandroff topology $(X, \mathcal{T})$, the {\em $T_0$-topogenous matrix A} \citep{topogenous} associated with it is defined as the $m \times n$ binary matrix defined as $A(i,j) = 1$ if $x_j \in U_i$, and $A(i,j) = 0$ otherwise. In words, each row defines the open sets that element $j$ belongs to, and each column defines the elements that are contained in open set $U_i$.
\end{definition}

\begin{theorem}
Given a finite space Alexandroff topology $(X, \mathcal{T})$, the {\em $T_0$-topogenous matrix A} associated with it defines a separating set. 
\end{theorem}

{\bf Proof:} Note that in a $T_0$ finite space Alexandroff topology, each element is distinguished by a unique neighborhood. Consequently, the open sets $U_x$ and $U_x$ associated with $x, y \in X$ must be distinct, and if $U_x = U_y$, then trivially $x = y$. Consequently, the topogenous matrix defines a separating set for the topology $(X, \mathcal{T})$. Similarly, the closed sets $F_x$ in Algorithm 1 also define separating sets.\qed

\subsection{Learning Causal Poset from Interventions} 

\begin{algorithm}[t]
\caption{Learn Causal Poset using Interventions and Conditional Independence Oracle.}
\SetAlgoLined
\KwIn{Dataset ${\cal D} = ({\cal E}, {\cal T})$ of a finite set of events ${\cal E}$ and their associated {\em open sets} $U \in {\cal T}$, represented as a topological space of open sets $g \subset {\cal E}$. A conditional independence oracle is also assumed.}
\KwOut{	Causally faithful poset ${\cal P}$ that is consistent with conditional independences in data.}
\Begin{
    Set the basic closed sets $F_e \leftarrow \emptyset$ \\ 
	\Repeat{convergence}{
		Select an open separating set $g \in {\cal T}$ and intervene on $g$. \\
	\For{$e \in g, f \notin g$}{
	 Use samples and the CI oracle to test if $(e \Perp f)_{{\cal M}_{g}}$ on dataset ${\cal D}$. \\ 
	 If CI test fails, then set $F_e \leftarrow F_e \cup \{ f \}$ because $f$ is an ancestor of $e$.
	 }
	 }
	Define the relation $e \leq f$ if $f \in F_e$, for all $e, f \in {\cal E}$, and compute its transitive closure.  \\ 
	Return the poset ${\cal P} = ({\cal E}, \leq)$, where $\leq$ is the induced relation on the poset ${\cal P}$. 
}
\end{algorithm}

\begin{algorithm}[h]
\caption{Given a causally faithful $T_0$ topological model $(X, \leq, {\cal T}$ and access to a conditional independence testing oracle, compute causal (observable) DAG model $G = (V,E)$.}
\SetAlgoLined
\KwIn{	 $T_0$  causally faithful topology ${\cal M} = (X, \leq, \cal{T})$  with $\leq$ defining a transitively closed partial ordering sets.}
\KwOut{	Causal DAG model $G = (V,E)$ on observable nodes. }
\Begin{
    Compute the {\em antichain} sets $T_i$ induced by partial ordering $\leq$. The set of open (or closed) sets in ${\cal M}$ are in $1-$ correspondence with the antichain sets.  \\ 
    Set $E \leftarrow \emptyset$. \\
	\Repeat{convergence}{
		Intervene on an antichain set $T_i$. \\
	\For{$x \in T_i, y \notin U$}{
	 Use samples and the CI oracle to test if $(x \Perp y)_{{\cal M}_{T_i}}$ on a dataset. \\ 
	 If CI test fails, then set $E \leftarrow E \cup (x,y)$.
	 }
	 }
    Set $V = X$, and return observable causal DAG $G = (V, E)$
}
\end{algorithm}

\begin{theorem}
Algorithm 1 requires only $| {\cal T}|$ interventions and conditional independence tests on samples obtained from each post-interventional distribution, to find a statistically consistent $T_0$ topological model. If there are $O(\log n)$ separating sets, the algorithm requires $O(\log n)$ interventions.
\end{theorem}

{\bf Proof:} If we intervene on the separating open  set $U$ and find an element $y \notin F_x$ that is statistically dependent on $x$, then we change $F_x$ to include $y$ (since $x$ needs to ``consult" $y$ in determining its value). The bound $O(\log n)$ in \citep{DBLP:conf/nips/KocaogluSB17} assumes that there are up to $2 \log n$ sets in the original separating system.  It has been argued in \citep{DBLP:conf/nips/BelloH18a} that interventions on multiple variables, such as used here and in the previous work \citep{DBLP:conf/nips/KocaogluSB17} can potentially require  an exponential number of experiments, if for example a separating set has $O(n)$ distinct elements, and each node is a binary variable, which requires two experiments (setting it to both $0$ and $1$). There is an inherent trade off between the size of each separating set, and the number of separating sets.  \qed 

Algorithm 2 is a generalization of Algorithm 1 in the recent paper by \citep{DBLP:conf/nips/KocaogluSB17}, who  construct a DAG by doing interventions on the {\em antichains} of posets, and extend this approach to discover causal models with latent variables as well. \citep{acharya} propose a related approach for inferring causal models, which does not require conditional independence testing, but uses a sample efficient testing methodology based on squared Hellinger distances.

Algorithm 2 constructs the observable DAG model, based on {\em antichain} sets, namely the set of incomparable elements at each level of the partial ordering. Antichain sets can be shown to be in bijective correspondence with the open sets on an Alexandroff $T_0$ topology.

\begin{theorem}
{\bf Mirsky's theorem} \citep{mirsky}: The {\bf height} of a $T_0$ topology causal model $(X, \leq, {\cal T})$ is defined to be the maximum cardinality of a chain, a totally ordered subset of the given partial order. For every partially ordered $T_0$ causal model $(X, \leq, {\cal T})$, the height also equals the minimum number of {\bf antichains}, namely subsets in which no pair of elements are ordered, into which the set may be partitioned.
\end{theorem}

\begin{theorem}
Algorithm $2$ requires $O(h)$ interventions and conditional independence tests on samples obtained from the post-interventional distributions, where $h$ is the height $h$ of a $T_0$ topology causal model $(X, \leq, {\cal T})$. 
\end{theorem}

Note that Algorithms 1 and 2 are topological generalizations of Algorithm 1 and 2 in  \citep{DBLP:conf/nips/KocaogluSB17}.   The remaining Algorithms 3 and 4 in \citep{DBLP:conf/nips/KocaogluSB17} on learning a latent variable DAG model can be generalized as well (see the appendix). 

\subsection{Efficient Enumeration of Homeomorphically Distinct Posets} 

Next, we turn to the fundamental problem of how to efficiently enumerate posets, which is a key requirement for scaling many causal discovery algorithms \citep{DBLP:conf/nips/KocaogluSB17,acharya,pmlr-v108-bernstein20a,BEERENWINKEL2006409}. 
\begin{definition}
For every $T_0$ finite space model ${\cal M}$  with a partial ordering $\leq$, define its associated {\bf Hasse diagram} $H_{\cal M}$ as a directed graph which captures all the relevant order information of ${\cal M}$. More precisely, the  vertices of $H_{\cal M}$ are the elements of ${\cal M}$, and the edges of $H_{\cal M}$ are such that there is a directed edge from $x$ to $y$ whenever $y \leq x$,  but there is no other vertex $z$ such that
$y \leq  z \leq x$.
\end{definition}

General pre-ordered finite spaces can be reduced to partially ordered $T_0$ topologies up to homomeomorphic equivalence. 

\begin{theorem}\citep{stong}
Let $(X, \cal{T})$ be an arbitrary finite space model with an associated preordering $\leq$. Let $X_0$ represent the quotient topological space $X / \sim$, where $x \sim y$ if $x \leq y$ and $y \leq x$. Then $X_0$ is a homotopically equivalent topological model with $T_0$ separability, and the quotient map $q: X \rightarrow X_0$ is a homotopy equivalence. Furthermore, $X_0$ induces a partial ordering on the elements $x \in X_0$. 
\end{theorem}

A key idea in the enumeration is to assume that each element in the Hasse diagram of the poset does not have an in-degree or out-degree of $1$. 
\begin{definition} \citep{stong}
An element $x \in X$ in a finite $T_0$  space $X$ is a {\bf down beat} point if $x$ covers one and only one element of of $X$. Alternatively, the set $\hat{U}_x = U_x \setminus \{x\}$ has a (unique) maximum. Similarly, $x \in X$ is an {\bf up beat} point if $x$ is covered by a unique element, or equivalently if $\hat{F}_x = F_x \setminus \{ x \}$ has a (unique) minimum. A {\bf beat} point is either a down beat or up beat point. 
\end{definition}
\begin{theorem}\citep{stong}
Let $X$ be a finite $T_0$ topological model, and let $x \in X$ be a (down, up) beat point. Then the reduced model $X \setminus \{ x \}$ is a {\bf strong deformation retract} of $X$.\footnote{In algebraic topology, a subspace $A \subset X$ is called a strong deformation retract of $X$ if there is a homotopy $F: X \times [0,1] \rightarrow A$ such that $F(x,0) = x, F(x, 1) \in A, F(a,t) = a$ for all $x \in X, t \in [0,1], a \in A$. \citep{munkres:algtop}.} A point x in a finite space ${\cal M}$ is an upbeat point if and only if it has
in-degree one in the associated Hasse diagram $H_{\cal M}$, i.e., it has only one incoming
edge). Similarly, $x$ is downbeat if and only if it has out-degree one (it has only one
outgoing edge).
\end{theorem} 

\begin{definition}
A finite $T_0$ topological space is a {\bf minimal model} if it has no beat points. A {\bf core} of a finite topological space $X$ is a strong deformation retract, which is a minimal finite space. The {\bf minimal graph} of a minimal model is its equivalent Hasse diagram. 
\end{definition}

\begin{theorem}\citep{stong}
{\bf Classification Theorem:} A homotopy equivalence between minimal finite space topological models is a homeomorphism. In particular, the core of a finite space model is unique up to homeomorphism and two finite spaces are homotopy equivalent if and only if they have homeomorphic cores. 
\end{theorem}

\begin{figure}
 \caption{Left: Constructing minimal posets by removing beat points \citep{barmak,stong}. Right: Efficiently enumerating minimal posets \citep{fix-patrias,brinkmann}.}
  \begin{minipage}[t]{.2\linewidth}
\vspace{0pt}
\centering
\includegraphics[scale=0.25]{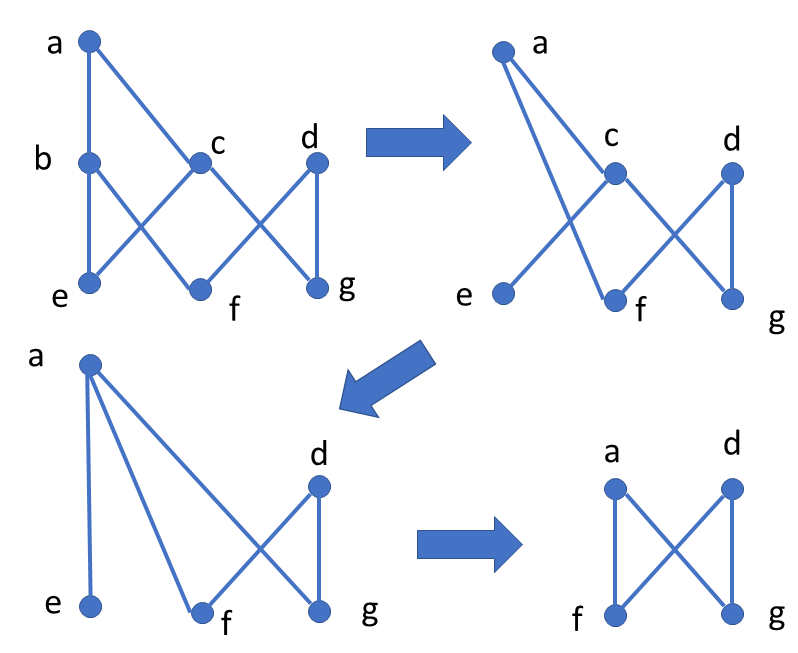} 
\end{minipage} \hfill 
 \begin{minipage}[t]{.5\linewidth}
\vspace{0pt}
\centering
\includegraphics[scale=0.13]{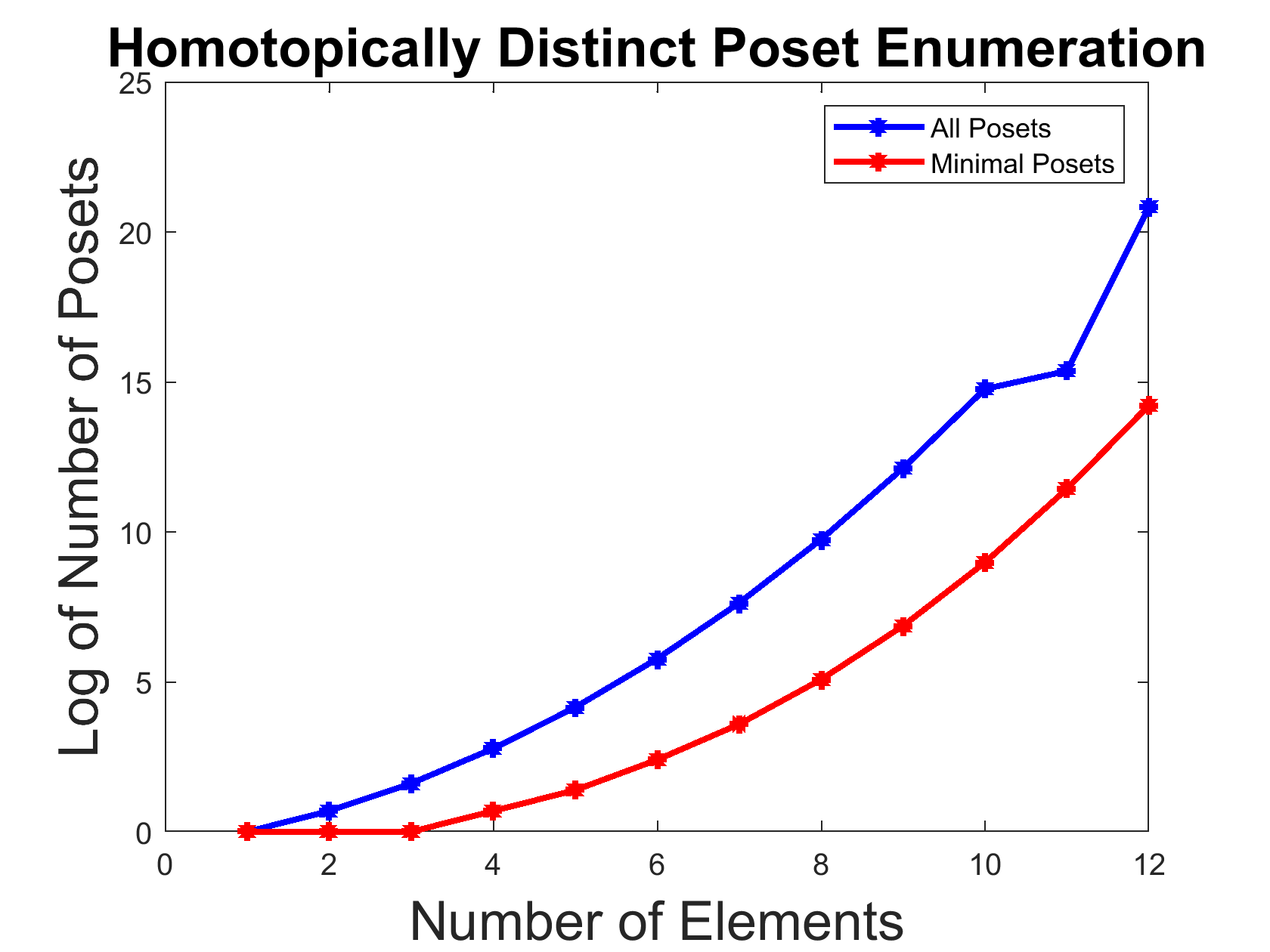} 
\end{minipage} 
 \label{enum-t0}
\end{figure}

Figure~\ref{enum-t0} illustrates the process of removing beat points to construct the minimal poset. b is an up beat point of $X$, c is an upbeat point of $X \setminus \{b \}$, and e is an up beat point of $X \setminus \{b, c\}$. Similarly, points c and e are removed, resulting in the minimal poset. The figure also shows that homeomorphic equivalences greatly reduces the search space of possible structures. Note the plot is on log scale.  For example, for $12$ variables, the number of minimal posets is $< 0.1$\% of the number of possible posets, a savings of three orders of magnitude. 

\begin{algorithm}[t]
\caption{Find Topologically Minimal $T_0$  Causal Model.}
\SetAlgoLined
\KwIn{	General pre-ordered causal model, such as a chain graph $G = (V,E)$ that is causally faithful to a dataset.\\
}
\KwOut{Minimal $T_0$ causal model homotopically equivalent to original non-$T_0$ model (e.g., from a chain graph $G$). \\
		The algorithm uses homotopy theory to find the {\em core} $T_0$ model of a general chain graph.}
\Begin{
    Define the topological model $(X, \mathcal{U})$ where $X = V$ and the open sets in $\mathcal{U}$ are constructed from the induced pre-order $\leq$ from $G$. 
	Define the minimal model $(X_0, \mathcal{U'})$, and set $X_0 = X$. \\
	\Repeat{convergence}{
	\For{$x,  y \in X_0$ s.t. $x \leq y, y \leq x$}{
	 Remove $x, y$ from $X_0$, and replace them with a new variable $z = x \sim y$. \\
	 Set $X_0 \leftarrow X_0 \setminus \{ x, y \} \cup \{ z \}$. $z$ represents the equivalence class that includes $x$ and $y$. \\
	 }
	 \For{$x \in X_0$}{
	 {\bf Remove down beat points:} If $\hat{U}_{x} = U_x \setminus \{ x \}$ has a maximum, then $X_0 \leftarrow X_0 \setminus \{ x \}$.\\ 
	 {\bf Remove up beat points:} If $\hat{F}_{x} = F \setminus \{ x \}$ has a minimum, then $X_0 \leftarrow X_0 \setminus \{ x \}$. 
	 }
	 }
	
	Define the open sets $U_x \in \mathcal{U'}$ as $U_x = \{y \ | \ y \leq x \}$ for $x \in X_0$. 
}
\end{algorithm}

Algorithm 3 determines a quotient $T_0$ topology that is homotopically equivalent to the general non-$T_0$ topology defined by a chain graphical models. Second, the algorithm further reduces the model to its core by removing {\em beat points} \citep{barmak, may,mccord,stong}. 
  
\section{Bioinformatics application}

 \begin{figure}[t]
\begin{minipage}[t]{.3\linewidth}
\vspace{0pt}
\centering
   \begin{tabular}{ll}
    Tumor & Gene   \\ \hline
    Pa017C & KRAS   \\ \hline 
    Pa017C & TP53   \\ \hline
    Pa019C & KRAS   \\ \hline
    Pa022C & KRAS   \\ \hline
    Pa022C & SMAD4   \\ \hline
    Pa022C & TP53  \\ \hline
    Pa032X & CDKN2A    \\ \hline
    \end{tabular}%
\end{minipage}%
\begin{minipage}[t]{.3\linewidth}
\vspace{0pt}
\centering
\includegraphics[scale=0.12]{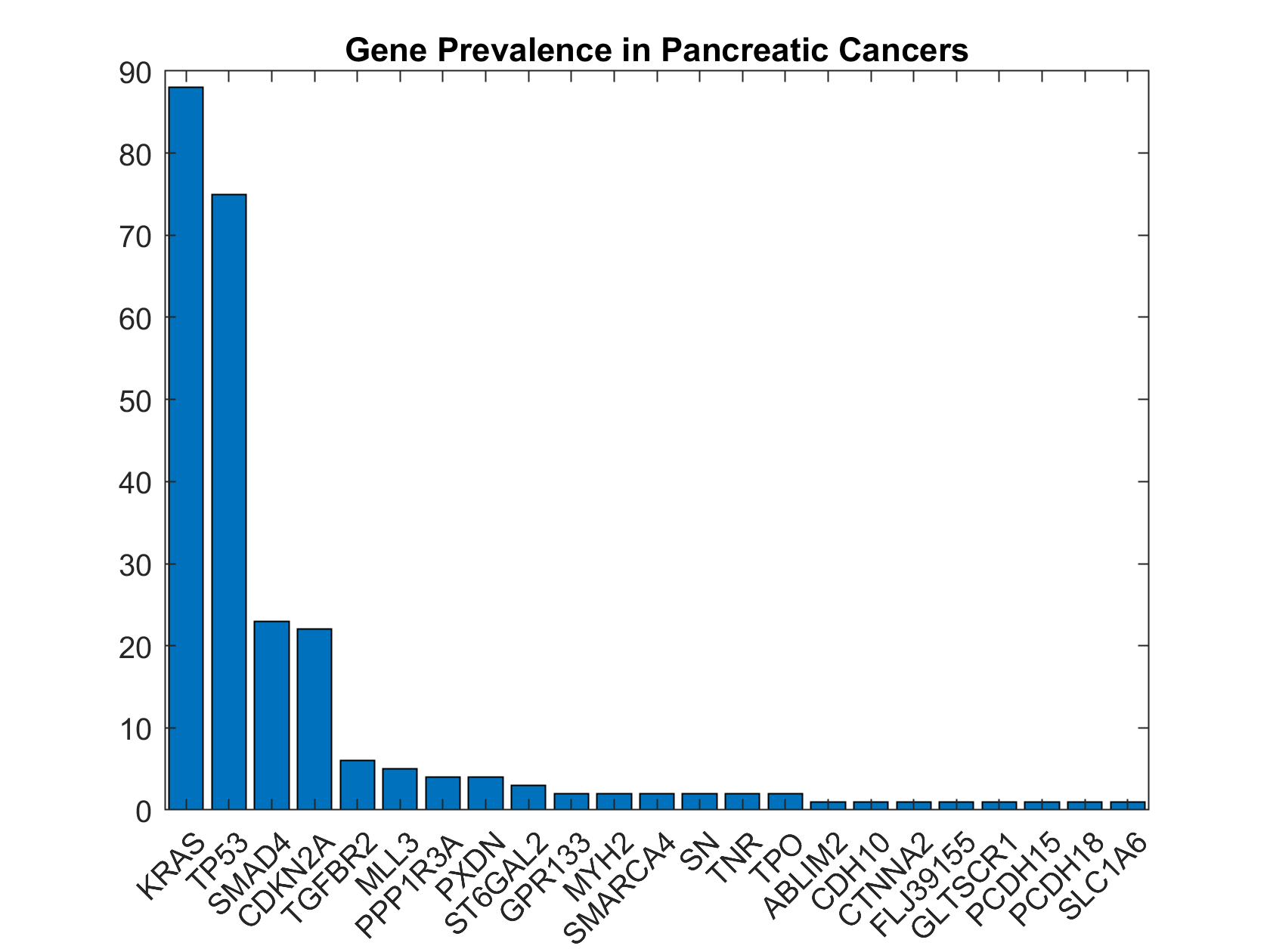} 
\end{minipage}
\begin{minipage}[t]{.3\linewidth}
\vspace{0pt}
\centering
\includegraphics[scale=0.12]{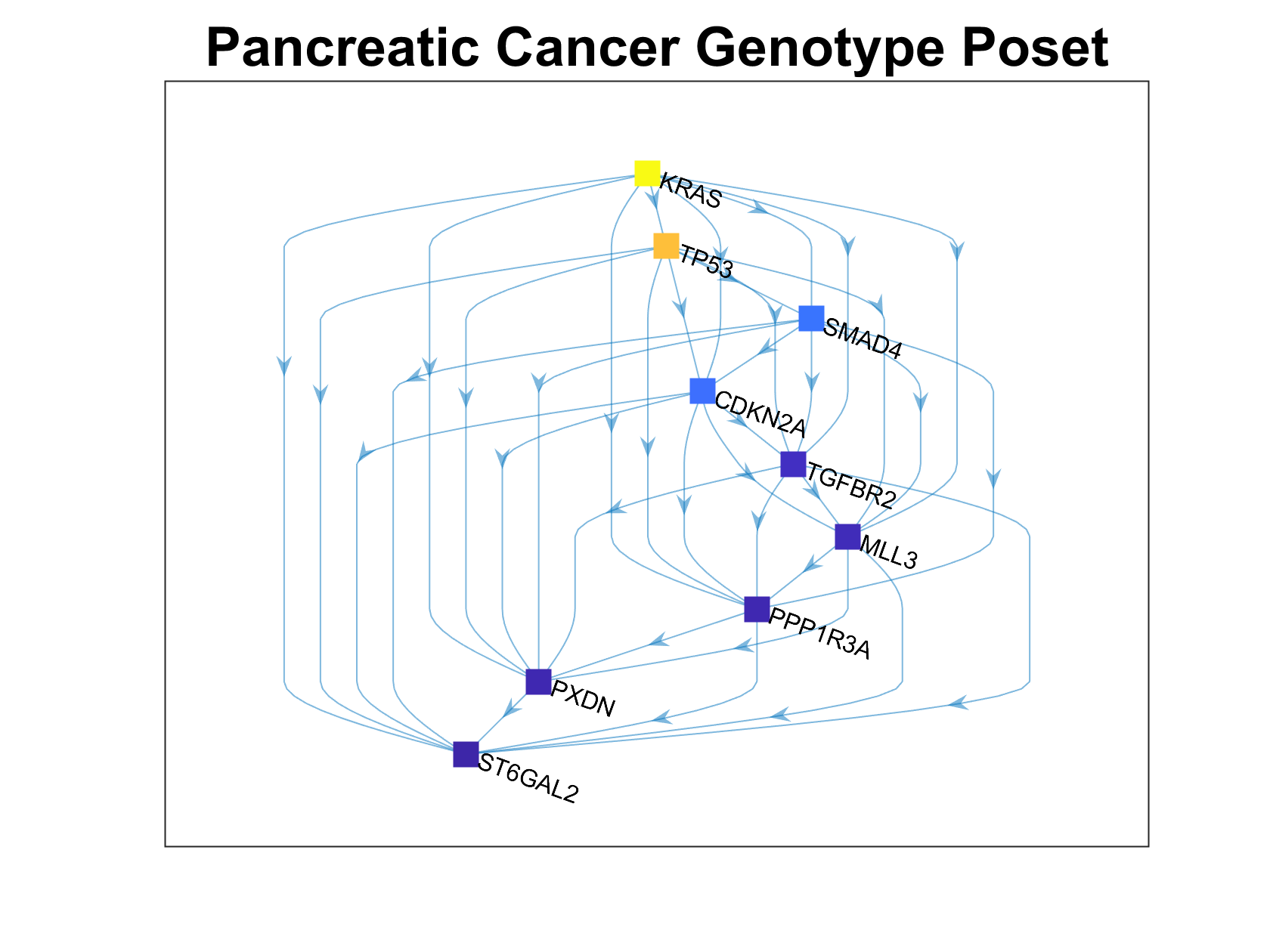} 
\end{minipage}
 \caption{Left: Genetic mutations in pancreatic cancer \citep{Jones1801}. Middle: histogram of genes sorted by mutation frequencies. Right: Poset learned from dataset.} 
\label{cancer-dataset}
\end{figure}%

Table~\ref{cancer-dataset} shows a small fragment of a dataset for pancreatic cancer \citep{Jones1801}. Like many cancers, it is marked by a particular partial ordering of mutations in some specific genes, such as {\bf KRAS}, {\bf TP53}, and so on.  In order to understand how to model and treat this deadly disease, it is crucial to understand the inherent partial ordering in the mutations of such genes. Pancreatic cancer remains one of the most prevalent and deadly forms of cancer. Roughly half a million humans contract the disease each year, most of whom succumb to it within a few years. \footnote{Sadly, this disease killed the much admired and long-time host of Jeopardy, Alex Trebek, last year.} Figure~\ref{cancer-dataset} shows the roughly $20$ most common genes that undergo mutations during the progression of this disease. The most common gene, the KRAS gene,  provides instructions for making a protein called K-Ras that is part of a signaling pathway known as the RAS/MAPK pathway. The protein relays signals from outside the cell to the cell's nucleus. The second most common mutation occurs in the TP53 gene, which  makes the p53 protein that normally acts as the supervisor in the cell as the body tries to repair damaged DNA. Like many cancers, pancreatic cancers occur as the normal reproductive machinery of the body is taken over by the cancer. 

In the pancreatic cancer problem, for example, the topological space $X$ is comprised of the significant events that mark the progression of the disease, as shown in Table~\ref{cancer-dataset}. In particular, the table shows that specific genes are mutated at specific locations by the change of an amino acid, causing the gene to malfunction. We can model a tumor in terms of its {\em genotype}, namely the subset of $X$, the gene events, that characterize the tumor. For example, the table shows the tumor {\bf Pa022C} can be characterized by the genotype {\bf KRAS}, {\bf SMAD4}, and {\bf TP53}. In general, a finite space topology is just the elements of the space (e.g. genetic events), and the subspaces (e.g., genomes) that define the topology.  

We illustrate our framework using the problem of inferring topological causal models for cancer  \citep{10.1093/biomet/asp023,cbn,BEERENWINKEL2006409,gerstung}. The progression of many types of cancer are marked by mutations of key genes whose normal reproductive machinery is subverted by the cancer \citep{Jones1801}. Often, viruses such as HIV and COVID-19 are constantly mutating to combat the pressure of interventions such as drugs, and successful treatment requires understanding the partial ordering of mutations. A number of past approaches use topological separability constraints on the data, assuming observed genotypes separate events, which as we will show, is abstractly a separability constraint on the underlying topological space.  

\begin{figure}[t]
\centering
\begin{minipage}{1\textwidth}
\centering
\includegraphics[scale=0.17]{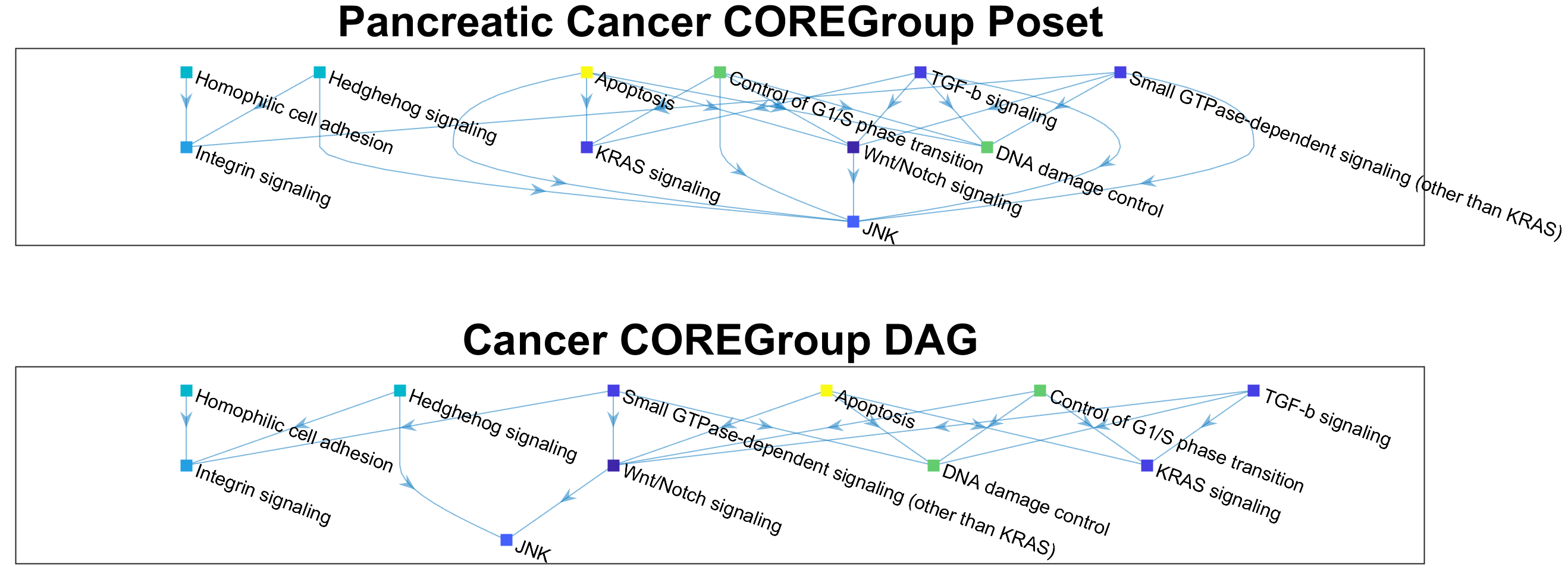} 
\includegraphics[scale=0.25]{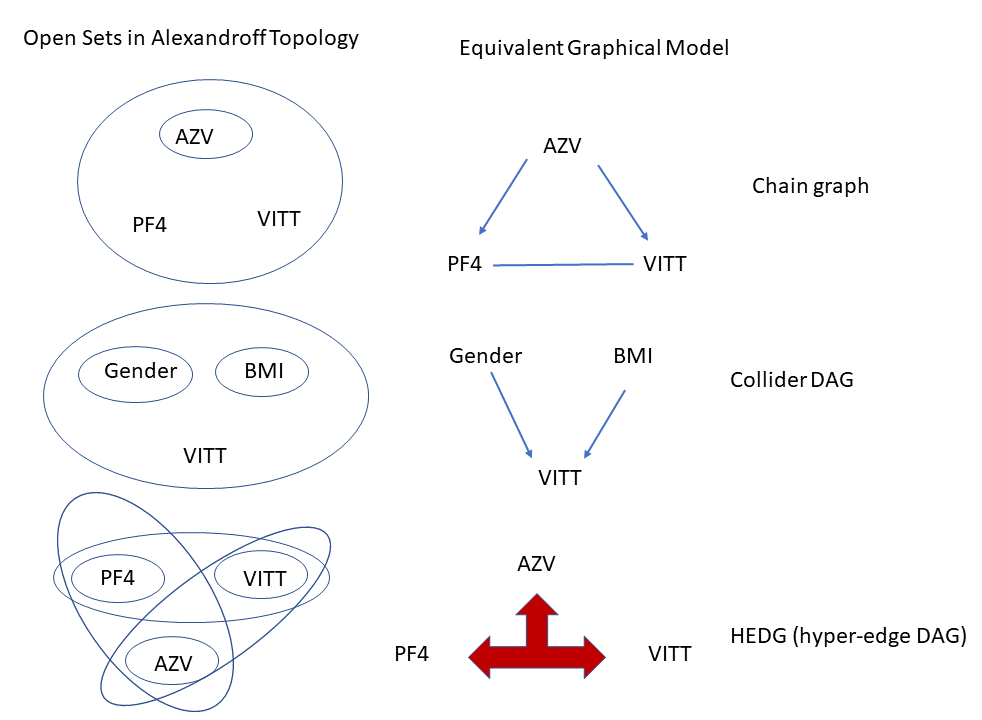}
\end{minipage}
\caption{Top: causal poset and DAG model of pathways in pancreatic cancer learned from a real-world dataset \citep{Jones1801}, showing genetic mutations occur along distinct pathways. Bottom: Topological representation of causal poset and DAG models for COVID-19  {\bf AZV} (AstraZeneca vaccine), {\bf VITT} (vaccine-induced clotting of blood), and other factors \citep{doi:10.1056/NEJMoa2104840,doi:10.1056/NEJMoa2104882}.}
\label{covid-diagram3}
\end{figure}

\begin{table}[htbp]
  \centering
  \caption{Core signaling pathways and processes genetically altered in most pancreatic cancers \citep{Jones1801}.}
    \begin{tabular}{|llll|} \hline
Regulatory pathway	& \% altered genes	& Tumors	& Representative altered genes \\ \hline
Apoptosis &	9	& 100\%	& CASP10, VCP, CAD, HIP1 \\ \hline 
DNA damage control	& 9 &	83\% &	ERCC4, ERCC6, EP300, TP53 \\ \hline
G1/S phase transition	& 19 &	100\% &	CDKN2A, FBXW7, CHD1, APC2 \\ \hline
Hedgehog signaling	& 19	& 100\%	& TBX5, SOX3, LRP2, GLI1, GLI3\\ \hline 
Homophilic cell adhesion &	30	& 79\%	& CDH1, CDH10, CDH2, CDH7, FAT\\ \hline 
Integrin signaling	& 24 & 	67\%	& ITGA4, ITGA9, ITGA11, LAMA1 \\ \hline
c-Jun N-terminal kinase &	9 & 96\%	& MAP4K3, TNF, ATF2, NFATC3 \\ \hline
KRAS signaling	& 5	& 100\%	& KRAS, MAP2K4, RASGRP3 \\ \hline
Regulation of invasion	& 46	& 92\%	& ADAM11, ADAM12, ADAM19\\ \hline 
Small GTPase–dependent & 33 &	79\%	& AGHGEF7, ARHGEF9, CDC42BPA\\ \hline 
TGF-$\beta$ signaling &	37	& 100\%	& TGFBR2, BMPR2, SMAD4, SMAD3 \\  \hline 
Wnt/Notch signaling	& 29 &	100\% &	MYC, PPP2R3A, WNT9A\\ \hline 
\end{tabular}
\label{pathways}
\end{table}

A key computational level in making model discovery tractable in evolutionary processes, such as pancreatic cancer, is that multiple sources of information are available that guide the discovery of the underlying poset model. In particular, for pancreatic cancer \citep{Jones1801}, in addition to the tumor genotype information show in Table~\ref{cancer-dataset}, it is also known that the disease follows certain pathways, as shown in Table~\ref{pathways}. This type of information from multiple sources gives the ability to construct multiple posets that reflect different event constraints \citep{BEERENWINKEL2006409}.

Algorithm 4 is a generalization of past algorithms that  infer conjunctive Bayesian networks (CBN)  from a dataset of events (e.g., tumors or signaling pathways) and their associated genotypes (e.g., sets of genes) \citep{cbn,BEERENWINKEL2006409}  The pathway poset and DAG shown in Figure~\ref{covid-diagram1} and the poset in Figure~\ref{cancer-dataset} were learned using Algorithm 4 using the pancreatic cancer dataset published in \citep{Jones1801}. 

\begin{algorithm}[t]
\caption{Application of Poset Discovery Algorithm to Bioinformatics.}
\SetAlgoLined
\KwIn{Dataset ${\cal D} = ({\cal E}, {\cal T})$ of a finite set of events ${\cal E}$ and their associated {\em genotypes} $U \in {\cal T}$, represented as a topological space of open sets $g \subset {\cal E}$. Here, it is assumed that each genome is an intervention target, whose size will affect the complexity of each causal experiment. A conditional independence oracle is also assumed.}
\KwOut{	Causally faithful poset ${\cal P}$ that is consistent with conditional independences in data.}
\Begin{
    Set the basic closed sets $F_e \leftarrow \emptyset$ \\ 
	\Repeat{convergence}{
		Select an open separating set $g \in {\cal T}$ and intervene on $g$. \\
	\For{$e \in g, f \notin g$}{
	 Use samples and the CI oracle to test if $(e \Perp f)_{{\cal M}_{g}}$ on dataset ${\cal D}$. \\ 
	 If CI test fails, then set $F_e \leftarrow F_e \cup \{ f \}$ because $f$ is an ancestor of $e$.
	 }
	 }
	Define the relation $e \leq f$ if $f \in F_e$, for all $e, f \in {\cal E}$, and compute its transitive closure.  \\ 
	Return the poset ${\cal P} = ({\cal E}, \leq)$, where $\leq$ is the induced relation on the poset ${\cal P}$. 
}
\end{algorithm}

\subsection{Limitations and Future Work}

We proposed a topological framework for causal discovery, building on the key relationship between posets and finite Alexandroff topologies. In the supplementary material, we elaborate on additional details.  We gave some examples from the domain of cancer genomics. A growing body of work in causal discovery has implicitly used topological constraints to constrain search. Our paper uses insights from algebraic topology of finite spaces into developing more scalable algorithms. Our paper has a number of significant limitations. We did not discuss building poset models over latent variables, which is important in many applications. Furthermore, a deeper study of the empirical performance of the algorithms proposed here is necessary to fully evaluate the promise of the proposed framework.


\end{document}